\documentclass{article}
\usepackage{amssymb,amsmath,amsfonts, latexsym}

\newtheorem{proposition}{Proposition}

\def\P{\ensuremath{\mathbb{P}}}

\def\R{\ensuremath{\mathbb{R}}}

\def\F{\ensuremath{\mathbb{F}}}

\def\int{\mathrm{int}}

\def\<{\ensuremath{\langle}}
\def\>{\ensuremath{\rangle}}

\begin{document}

\title{Topology of real toric surfaces}

\author{Sam Payne}

\maketitle

\begin{abstract}
    
We determine the homeomorphism type of the set of real points of a
smooth projective toric surface.\footnote{When this note was prepared
and submitted to the arXiv, the author was not aware that these
results (and much more) had already appeared in C. Delaunay's work on
real toric varieties \cite{Delaunay} \cite{DelThesis}.  We hope that
this note may serve as an expository introduction to some of the ideas
and techniques in Delaunay's work.}
\end{abstract}

The purpose of this note is to prove the following topological 
classification of real toric surfaces.

\vspace{10 pt}

    \noindent \textbf{Theorem} \emph{\cite[Proposition 
    4.1.2.]{Delaunay} Let $X = X(\Delta)$ be a smooth
    projective toric surface.  If $X$ is isomorphic to an even
    Hirzebruch surface $\F_{2a}$, for some integer $a \geq 0$, then
    $X(\R)$ is homeomorphic to $S^{1} \times S^{1}$.  Otherwise,
    $X(\R)$ is homeomorphic to a connect sum of $\# \Delta(1) - 2$
    copies of $\R\P^{2}$.}

\vspace{10 pt}

\noindent Here, $\Delta(1)$ denotes the set of 1-dimensional cones in
$\Delta$.

We begin with a few preliminaries on the topology of real projective
toric varieties in general, roughly following \cite[Chapter 11,
Section 5]{GKZ}.  Let $X = X(\Delta)$ be a projective toric variety,
and let $P \subset M_{\R}$ be the polytope corresponding to an ample
toric divisor on $X$.  For each group homomorphism $\epsilon : M
\rightarrow \{ \pm 1 \}$, define
\[
      T_{\epsilon} := \{ t \in T(\R) \subset X(\R) : \mathrm{sgn} 
      (\chi^{u}(t)) = \epsilon(u) \mbox{ for all } u \in M \}.
\]      
Let $X_{\epsilon} \subset X(\R)$ be the closure of $T_{\epsilon}$ in
the real analytic topology.  Let 
\[ 
      \mu : X(\R) \rightarrow P
\]
be the moment map, and $\mu_{\epsilon}$ the restriction of $\mu$ to
$X_{\epsilon}$.  We claim that $\mu_{\epsilon}$ is a homeomorphism. 
The case $\epsilon_{0}(M) = +1$ is proved in \cite[Section
4.2]{Fulton}.  For general $\epsilon,$ the semigroup homomorphism
$\epsilon: M \rightarrow \R$ corresponds to a point $t_{\epsilon} \in
T(\R)$.  Translation by $t_{\epsilon}$ takes $X_{\epsilon_{0}}$ to
$X_{\epsilon}$ and commutes with the moment map, so the claim follows.

Now $X(\R) = \bigcup_{\epsilon} X_{\epsilon}$.  The following
proposition describes $X_{\epsilon} \cap X_{\epsilon'}$ and the
induced construction of $X(\R)$ by gluing $2^{\dim X}$ copies
$P_{\epsilon}$ of $P$, indexed by the group homomorphisms $\epsilon: 
M \rightarrow \{ \pm 1 \}$.  For a face $F \subset P$, let $F_{\epsilon}$ 
denote the corresponding face of $P_{\epsilon}$, and let $M_{F}$ be the
subgroup of $M$ parallel to $F$, i.e. \ if $\tau \in \Delta$ is the cone
corresponding to $F$, then $M_{F} = \tau^{\perp} \cap M$.

\begin{proposition} \label{gluing}
    \cite[Proposition 4.1.1.]{Delaunay} Let $X = X_{P}$ be a
    projective toric variety, and let $\epsilon, \epsilon': M
    \rightarrow \{ \pm 1 \}$ be group homomorphisms.  Then
    $X_{\epsilon} \cap X_{\epsilon'}$ is the union of the preimages
    under $\mu_{\epsilon}$ of the faces $F \subset P$ such that
    $\epsilon|_{M_{F}} = \epsilon'|_{M_{F}}$.  In particular, $X(\R)$
    is homeomorphic to the space constructed by gluing the
    $P_{\epsilon}$ along faces as follows: $F_{\epsilon}$ and
    $F_{\epsilon'}$ are identified if and only if $\epsilon|_{M_{F}} =
    \epsilon'|_{M_{F}}$.
\end{proposition}

\noindent Proof: Let $F$ be a face of $P$, and let $F^{\circ}$ be the
relative interior of $F$.  Let $x \in \mu_{\epsilon}^{-1}(F^{\circ})$. 
It will suffice to show that $x \in X_{\epsilon'}$ if and only if
$\epsilon|_{M_{F}} = \epsilon'|_{M_{F}}$.  The rational functions $\{
\chi^{u}: u \in M_{F} \}$ are regular on $\mu^{-1}(F^{\circ})$ and
separate points.  For $u \in M_{F}$, the absolute value of
$\chi^{u}(x)$ is determined by $\mu(x)$ and the sign of $\chi^{u}(x)$
is $\epsilon(u)$.  Hence $x \in X_{\epsilon'}$ if and only if
$\epsilon'(u) = \epsilon(u)$ for all $u \in M_{F}$.  \nolinebreak
\hfill $\Box$

\vspace{10 pt}

Proposition \ref{gluing} is a correction of Theorem 5.4 from
\cite[Chapter 11]{GKZ}, which says that $F_{\epsilon}$ and
$F_{\epsilon'}$ are identified if and only if $\epsilon$ and
$\epsilon'$ agree on the intersection of $M$ with the affine span of
$F$.  If $\epsilon$ and $\epsilon'$ agree on $M \cap \mathrm{Aff}(F)$,
then they also agree on $M_{F}$, since every point in $M_{F}$ is a
difference of points in $M \cap \mathrm{Aff}(F)$, but the converse is
false in general.

The gluing construction does not depend on the choice of ample
divisor.  Indeed, $P$ can be replaced by the unit ball in $N_{\R}$
with a cell structure on the boundary sphere dual to that induced by
intersection with the cones of $\Delta$.  If $F$ is the cell
corresponding to a cone $\tau \in \Delta$, then $F_{\epsilon}$ and
$F_{\epsilon'}$ are identified if and only if $\epsilon$ and
$\epsilon'$ agree on $\tau^{\perp} \cap M$.

\vspace{10 pt}

\noindent \textbf{Example} \ Suppose $X \cong \F_{a}$ is a Hirzebruch
surface, where $a$ is a nonnegative integer.  Choose coordinates on
$N$ so that the rays of $\Delta$ are generated by $e_{1}, e_{2},
-e_{2}$, and $-e_{1} + a e_{2}$.  To visualize the construction of
$X(\R)$ by the gluing recipe in Proposition \ref{gluing}, draw four
copies $P_{\epsilon}$ of the polytope $P$ corresponding to some ample
divisor on $X$, and label the interior of $P_{\epsilon}$ with an
ordered pair of plus or minus signs for $\epsilon(e_{1}^{*}),
\epsilon(e_{2}^{*})$.  Similarly, label each edge $F_{\epsilon}$ of
$P_{\epsilon}$ with a plus or minus sign for $\epsilon(u_{F})$, where
$u_{F}$ is a primitive generator of the 1-dimensional lattice 
$M_{F}$.  The resulting four labeled copies of $P$ are as follows:

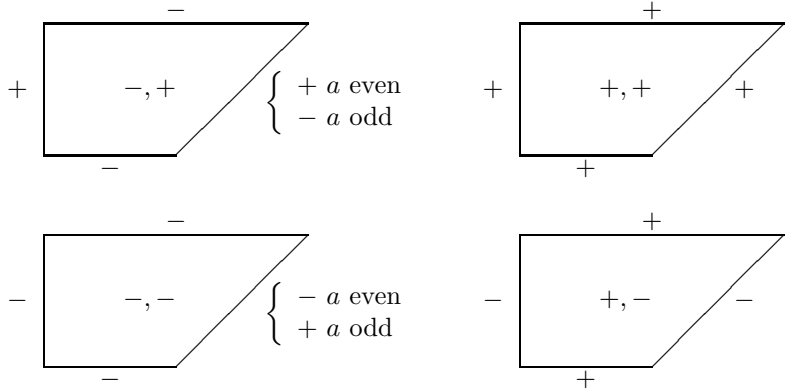
\begin{picture}(240, 150)(-30,0)
    \put(0,0){\line(1,0){50}}
    \put(0,0){\line(0,1){50}}
    \put(50,0){\line(1,1){50}}
    \put(0,50){\line(1,0){100}}
    
    \put(-10,25){\makebox(0,0){$-$}} 
    \put(25,-5){\makebox(0,0){$-$}}
    \put(80,20){\makebox(0,0)[l]{ $\left\{ \begin{array}{ll} - 
    \mbox{  $a$
               even} \\ + \mbox{  $a$ odd} \end{array} \right.$} }
    \put(50,55){\makebox(0,0){$-$}}
    
    \put(40,25){\makebox(0,0){$-,-$}}
    
    \put(180,0){\line(1,0){50}}
    \put(180,0){\line(0,1){50}}
    \put(230,0){\line(1,1){50}}
    \put(180,50){\line(1,0){100}}
    
    \put(170,25){\makebox(0,0){$-$}}
    \put(205,-5){\makebox(0,0){$+$}}
    \put(265,25){\makebox(0,0){$-$}}
    \put(230,55){\makebox(0,0){$+$}}
    
    \put(220,25){\makebox(0,0){$+,-$}}

    \put(180,80){\line(1,0){50}}
    \put(180,80){\line(0,1){50}}
    \put(230,80){\line(1,1){50}}
    \put(180,130){\line(1,0){100}}
    
    \put(170,105){\makebox(0,0){$+$}}
    \put(205,75){\makebox(0,0){$+$}}
    \put(265,105){\makebox(0,0){$+$}}
    \put(230,135){\makebox(0,0){$+$}}
    
    \put(220,105){\makebox(0,0){$+,+$}}

    \put(0,80){\line(1,0){50}}
    \put(0,80){\line(0,1){50}}
    \put(50,80){\line(1,1){50}}
    \put(0,130){\line(1,0){100}}
    
    \put(-10,105){\makebox(0,0){$+$}} 
    \put(25,75){\makebox(0,0){$-$}}
    \put(80,100){\makebox(0,0)[l]{ $\left\{ \begin{array}{ll} + 
    \mbox{  $a$
               even} \\ - \mbox{  $a$ odd} \end{array} \right.$} }
    \put(50,135){\makebox(0,0){$-$}}
    
    \put(40,105){\makebox(0,0){$-,+$}}

\end{picture}

\noindent The upper left copy of $P$ represents
$P_{\epsilon}$, where $\epsilon(e_{1}^{*}) = -1
$ and $\epsilon(e_{2}^{*}) =
+1$.  The right (diagonal) edge of $P_{\epsilon}$ is marked with + if 
$a$ is even and with $-$ if $a$ is odd because it is parallel to
the subgroup generated by $ae_{1}^{*} + e_{2}^{*}$, and
\[
    \epsilon(ae_{1}^{*} + e_{2}^{*}) \ = \ \epsilon(e_{1}^{*})^{a} \cdot 
    \epsilon(e_{2}^{*}) \ = \ (-1)^{a}.
\]    
With this notation, $F_{\epsilon}$ and $F_{\epsilon'}$ are identified
if and only if they are marked with the same sign.  Gluing the top,
bottom, and left edges produces a cylinder whose ends are formed by
the right edges.  The way the ends of the cylinder are glued depends
on the parity of $a$.  If $a$ is even the result is a torus; if $a$ is
odd, the result is a Klein bottle.

\vspace{10 pt}

Let $\Delta(i)$ be the set of $i$-dimensional cones of $\Delta$.

\begin{proposition} \label{Euler characteristic}
    The Euler characteristic of $X(\R)$ is given by
    \[
       \chi(X(\R)) \ = \ \sum_{i=0}^{\dim X} (-2)^{\dim X - i} \cdot \# 
       \Delta(i).
    \]
\end{proposition}

\noindent Proof: Proposition \ref{gluing} gives a cell decomposition
of $X(\R)$ with $2^{\dim X - i} \cdot \# \Delta(i)$ cells of dimension
$\dim X - i$.  \hfill $\Box$
    
\vspace{10 pt}

For the remainder of this note, we assume that $X$ is a smooth 
projective toric surface.

\vspace{10 pt}

\noindent \textbf{Lemma} \emph{
    Let $C$ be a $T$-invariant curve in $X$, and let $U$ be a tubular
    neighborhood of $C(\R)$ in $X(\R)$.  If $C^{2}$ is even then $U$
    is homeomorphic to a cylinder; if $C^{2}$ is odd, then $U$ is
    homeomorphic to a M\"{o}bius band.}

\vspace{10 pt}

\noindent Proof: Say $C = V(\rho)$ and $C^{2} = a$.  We can choose
coordinates on $N$ so that $\rho$ is generated by $e_{2}$ and the
adjacent rays in $\Delta$ are generated by $e_{1}$ and $-e_{1} - a
e_{2}$.  Then $\mu(C)$ is the bottom edge of $P$ (the edge whose
inward normal is generated by $e_{2}$) and a tubular neighborhood
of $C(\R)$ may be constructed by gluing neighborhoods of the bottom
edges of the $P_{\epsilon}$ according to the recipe given in
Proposition \ref{gluing}.  With notation as in the Example, the pieces
to be glued are drawn in the figure below.

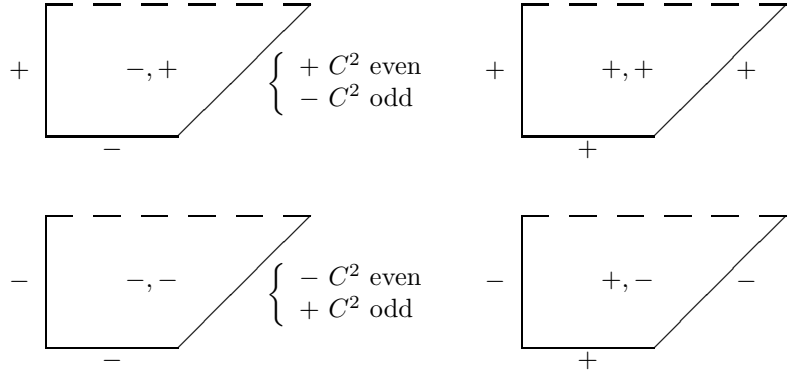
\begin{picture}(240, 140)(-30,0)
    \put(0,0){\line(1,0){50}}
    \put(0,0){\line(0,1){50}}
    \put(50,0){\line(1,1){50}}
    \multiput(0,50)(18,0){6}{\line(1,0){10}}
    
    \put(-10,25){\makebox(0,0){$-$}} 
    \put(25,-5){\makebox(0,0){$-$}}
    \put(80,20){\makebox(0,0)[l]{ $\left\{ \begin{array}{ll} - 
    \mbox{  $C^{2}$
               even} \\ + \mbox{  $C^{2}$ odd} \end{array} \right.$} }
    
    \put(40,25){\makebox(0,0){$-,-$}}
    
    \put(180,0){\line(1,0){50}}
    \put(180,0){\line(0,1){50}}
    \put(230,0){\line(1,1){50}}
    \multiput(180,50)(18,0){6}{\line(1,0){10}}
    
    \put(170,25){\makebox(0,0){$-$}}
    \put(205,-5){\makebox(0,0){$+$}}
    \put(265,25){\makebox(0,0){$-$}}
    
    \put(220,25){\makebox(0,0){$+,-$}}

    \put(180,80){\line(1,0){50}}
    \put(180,80){\line(0,1){50}}
    \put(230,80){\line(1,1){50}}
    \multiput(180,130)(18,0){6}{\line(1,0){10}}
    
    \put(170,105){\makebox(0,0){$+$}}
    \put(205,75){\makebox(0,0){$+$}}
    \put(265,105){\makebox(0,0){$+$}}
    
    \put(220,105){\makebox(0,0){$+,+$}}

    \put(0,80){\line(1,0){50}}
    \put(0,80){\line(0,1){50}}
    \put(50,80){\line(1,1){50}}
    \multiput(0,130)(18,0){6}{\line(1,0){10}}
    
    \put(-10,105){\makebox(0,0){$+$}} 
    \put(25,75){\makebox(0,0){$-$}}
    \put(80,100){\makebox(0,0)[l]{ $\left\{ \begin{array}{ll} + 
    \mbox{  $C^{2}$
               even} \\ - \mbox{  $C^{2}$ odd} \end{array} \right.$} }
    
    \put(40,105){\makebox(0,0){$-,+$}}

\end{picture}

\noindent Gluing the bottom and left edges yields a strip homeomorphic
to $[0,1] \times (0,1)$ whose ends are formed by the right edges.  The
gluing of the ends of the strip depends on the parity of $a$.  If $a$
is even the result is a cylinder; if $a$ is odd the result is a
M\"{o}bius band.  \hfill $\Box$

\vspace{10 pt}

\noindent Proof of Theorem: The homeomorphism type of a surface is
determined by its Euler characteristic and orientability.  By
Proposition \ref{Euler characteristic}, $\chi(X(\R)) = 4 - \#
\Delta(1)$. It remains to show that $X(\R)$ is orientable if and
only if $X$ is isomorphic to an even Hirzebruch surface $\F_{2a}$.

If $X$ is not a minimal surface, then it contains a $-1$-curve $C$,
which must be $T$-invariant (otherwise $C$ would move and hence have
nonnegative self-inter\-section).  By the Lemma, a tubular
neighborhood of $C(\R)$ is homeomorphic to a M\"{o}bius band, so
$X(\R)$ is nonorientable.  Therefore, if $X(\R)$ is orientable, then
$X$ must be minimal.  The minimal rational surfaces are $\P^{2}$ and
the Hirzebruch surfaces.  Of course, $\P^{2}(\R)$ is homeomorphic to
$\R\P^{2}$, which is nonorientable.  As seen in the Example,
$\F_{a}(\R)$ is orientable if and only if $a$ is even.  \hfill $\Box$

\vspace{10 pt}

\noindent Acknowledgments -- This note grew out of examples discussed
during an evening workshop on real toric varieties at the 2004 PCMI
conference on geometric combinatorics.  I wish to thank my fellow
participants J. Martin and V. Cormani, as well as B. Sturmfels, who
suggested the topic and led the discussion.

\end{document}